\documentclass[10pt]{article}
\usepackage{latexsym,amssymb,amsmath}

\def\N{\mathbb{N}}
\def\Q{\mathbb{Q}}

\def\R{\mathbb{R}}

\def\sus{\subset}
\def\al{\alpha}

\def\ep{\varepsilon}

\def\ds{\dots}

\begin{document}
\title{Birch's theorem: if $f(n)$ is multiplicative and has a non-decreasing normal order then $f(n)=n^{\al}$}
\author{Martin Klazar\footnote{{\tt klazar@kam.mff.cuni.cz}}}

\maketitle

\begin{abstract}
For pedagogical purposes (inclusion in lecture notes) we review the proof of the theorem
stated in the title. At the end we state a problem.
\end{abstract}

\section{Introduction}

In 1967 B.\,J. Birch, later of the Birch and Swinnerton-Dyer conjecture fame, proved in \cite{birc} a most interesting result.

\bigskip\noindent
{\bf Theorem (Birch, 1967). }{\em The only multiplicative functions $f:\N\to\R_{\ge0}$ that are unbounded and have a non-decreasing 
normal order are the powers of $n$, the functions $f(n)=n^{\al}$ for a constant $\al>0$.  
}

\bigskip\noindent
Multiplicativity means that $f(mn)=f(m)f(n)$ for every two coprime numbers $m,n\in\N$ (thus $f(1)=1$ unless $f\equiv0$), 
$\N=\{1,2,\ds\}$, and the clause about a non-decreasing normal order means that a non-decreasing function 
$g:\N\to\R_{>0}$ exists such that for every $\ep>0$, $\#(n\le x\;|\;\frac{f(n)}{g(n)}\not\in(1-\ep,1+\ep))=o(x)$ as $x\to+\infty$.

In this write-up I present the proof of Birch's theorem, as given in Birch \cite{birc} and Narkiewicz \cite[pp. 98--102]{nark} 
(see also \cite{nark_e}). It is a beautiful proof in the erd\H osian style. To be honest, I started with the intention to correct two errors 
I thought I had discovered in the argument. Fortunately, in the process of writing everything clarified and the errors disappeared.
Still, I will point out the two steps I struggled with. To the interested reader, much smarter than me, they will certainly pose 
no difficulty. 

\section{The proof with two conundrums}

We use notation of \cite{birc}, so let 
$$
b(n)=\log f(n)\ \mbox{ and }\ c(n)=\log g(n)\;. 
$$
Birch \cite[p. 149]{birc} writes just ``If $f$ is 
unbounded, then $g(n)$ tends to infinity with $n$, so we may suppose that $c(n)>0$ for all $n$.'' but Narkiewicz 
\cite[Lemat~2.5 on p. 98]{nark} gives more details. Assume for contrary that $g(n)$ has a finite limit $a>0$. Then, by the relation 
bounding $f$ and $g$, there are constants $0<A<a<B$ such that for every $x>0$ and $n\le x$ we have $A<f(n)<B$, with $o(x)$ exceptions.
Let $E\sus\N$ be the exceptions; $E$ has density $0$. Fix any $M>B$. Since $f$ is unbounded, there is an $m\in\N$ with $f(m)>M/A$. 
The sets $\{nm+1\;|\;n\in\N\}$ and $\{(nm+1)m\;|\;n\in\N\}$ have positive densities and thus so has 
$X=\{n\in\N\;|\;nm+1,(nm+1)m\not\in E\}$. For any $n\in X$ we get the contradiction $B>f((nm+1)m)=f(nm+1)f(m)>Af(m)>M$. 

Thus indeed $\lim g(n)=+\infty$. Changing finitely many values of $g(n)$ we may assume that always $g(n)>1$ and $c(n)>0$. 
By Birch \cite{birc}, ``Using the three conditions
\begin{quote}
given $\ep>0$, $|b(n)-c(n)|<\ep$ for all but $o(x)$ integers $n<x$;\\
$b(mn)=b(m)+b(n)$ if $(m,n)=1$;\\
$c(n)\ge c(m)>0$ for $n\ge m$;
\end{quote}
we gradually deduce more and more till everything collapses.'' Let $m,n\in\N$ and $\ep>0$ be arbitrary with 
$|b(m)-c(m)|,|b(n)-c(n)|<\ep$. We assume that $m,n\ge2$. 
It follows that for any $\eta\in(0,\frac{1}{2})$ there is an $S>0$ such that for every $R\ge S$ there are 
$s,t\in\N$ satisfying
$$
(1-\eta)R<s<R<t<(1+\eta)R,\ s\equiv t\equiv 1\ (\mathrm{mod}\ mn)
$$
and 
$$
|b(s)-c(s)|,\;|b(ms)-c(ms)|,\;|b(t)-c(t)|,\;|b(nt)-c(nt)|<\ep\;.
$$
(Only $o(R)$ of the integers $s\in((1-\eta)R,R)$ violate the first or the second lastly displayed inequality, and so
for large $R$ we certainly find there an $s\equiv1\ (\mathrm{mod}\ mn)$ satisfying both. The same for $t$.) From 
$b(ms)=b(m)+b(s)$ and $b(nt)=b(n)+b(t)$ we get
$$
|c(ms)-c(m)-c(s)|,\;|c(nt)-c(n)-c(t)|<3\ep\;.
$$

We define by induction numbers $s_0<s_1<\ds$ and $t_0<t_1<\ds$ in $\N$, all congruent to $1$ modulo $mn$, such that 
$$
(1-\eta)S<s_0<S<t_0<(1+\eta)S
$$ 
and, for every $i,j\in\N_0$,
$$
(1-\eta)ms_i<s_{i+1}<ms_i,\ nt_j<t_{j+1}<(1+\eta)nt_j\;,
$$
and
$$
|b(s_i)-c(s_i)|,\;|b(ms_i)-c(ms_i)|,\;|b(t_j)-c(t_j)|,\;|b(nt_j)-c(nt_j)|<\ep\;.
$$
(In the previous claim we first set $R=S$ and get $s_0=s$, then we set $R=ms_0(\ge S)$ and get $s_1=s$, and so on. 
Since $m\ge2$ and $\eta<\frac{1}{2}$, we stay above $S$ and $s_i$ increase. Similarly and more easily for $t_j$.) Then, as we know, 
for every $i\in\N_0$ one has
$$
|c(ms_i)-c(m)-c(s_i)|<3\ep\;.
$$
Monotonicity of $c$ gives 
$$
c(s_i)>c(ms_i)-c(m)-3\ep\ge c(s_{i+1})-c(m)-3\ep
$$
and so $c(s_h)<c(S)+hc(m)+3h\ep$ for every $h\in\N$ by iteration. On the other hand, $s_h>(1-\eta)^{h+1}m^hS$ by iterating the 
above inequalities. Similarly for $t_j$ we get $c(t_k)>c(S)+kc(n)-3k\ep$ for every $k\in\N$ and $t_k<(1+\eta)^{k+1}n^kS$. 

Now if $h,k\in\N$ are such that $m^h>n^k$, equivalently $h\log m>k\log n$ (recall that $\log m\ne0$), we may select $\eta>0$
so small that still 
$$
(1-\eta)^{h+1}m^h>(1+\eta)^{k+1}n^k\;.
$$
This implies that $s_h>t_k$ and $c(s_h)\ge c(t_k)$ (by monotonicity of $c$), hence $hc(m)+3h\ep>kc(n)-3k\ep$ and
$$
\frac{h}{k}>\frac{c(n)-3\ep}{c(m)+3\ep}\;.
$$
It follows that 
$$
\frac{\log n}{\log m}\ge\frac{c(n)-3\ep}{c(m)+3\ep}\;.
$$
(But how come? {\em This is the first step I struggled with.} Don't we assume that $h/k>(\log n)/(\log m)$? To combine inequalities 
by transitivity we would need this one be opposite!)

Nevertheless, we get
$$
\frac{c(n)}{\log n}-\frac{c(m)}{\log m}\le3\ep\left(\frac{1}{\log m}+\frac{1}{\log n}\right)
$$
and, changing the roles of $m$ and $n$, the reverse inequality $\ds\ge-3\ep\ds\;$. So we have proved that
$$
\left|\frac{c(n)}{\log n}-\frac{c(m)}{\log m}\right|\le3\ep\left(\frac{1}{\log m}+\frac{1}{\log n}\right)
$$
whenever $|b(m)-c(m)|<\ep$ and $|b(n)-c(n)|<\ep$. This implies 
$$
\left|\frac{c(n)}{\log n}-\frac{c(m)}{\log m}\right|\le(|b(m)-c(m)|+|b(n)-c(n)|)\left(\frac{3}{\log m}+\frac{3}{\log n}\right)
$$
for all $m,n$. (But how come? {\em This is the second step I struggled with.} Let's say that the penultimate displayed 
inequality holds for every $m,n$ as an equality for $3\ep$ replaced with $2\ep$, and that we have $m,n$ such that 
$|b(m)-c(m)|,|b(n)-c(n)|<\ep/4$. The last two displayed inequalities then contradict each other!). 

Nevertheless, we conclude the proof. Obviously, $|b(n_i)-c(n_i)|\to0$ for a sequence $n_1<n_2<\ds\;$. The last displayed inequality shows that 
the values $c(n_i)/\log n_i$ are bounded. Passing to a subsequence we get $\lim_ic(n_i)/\log n_i=\al$, with a finite limit $\al$.
Setting $n=n_i$ and letting $i\to\infty$ gives
$$
|c(m)-\al\log m|\le3|b(m)-c(m)|\ \mbox{ and }\ |b(m)-\al\log m|\le4|b(m)-c(m)|
$$
for every $m\in\N$ (well, $m\ge2$). Thus, given any $\ep>0$, $|b(m)-\al\log m|<\ep$ for all but $o(x)$ numbers $m\le x$. Let
$E\sus\N$ be the set of exceptional $m$; it has density  $0$. We take any $m\in\N$. The set $X=\{n\in\N\;|\;(n,m)=1,n,mn\not\in E\}$ has positive
density. For any $n\in X$ we have
$$
|b(n)-\al\log n|,\;|b(mn)-\al\log(mn)|<\ep\;.
$$
So, by the additivity of the functions $b$ and $\log$, $\ep>|b(mn)-\al\log(mn)|\ge|b(m)-\al\log m|-|b(n)-\al\log n|$ and 
$|b(m)-\al\log m|<2\ep$. As this holds for any $\ep>0$, we get the desired equality 
$$
b(m)=\al\log m\ \mbox{ or }\ f(m)=m^{\alpha}
$$
for every $m\in\N$. We are done. Well, $\ds$

\section{Concluding remarks}

How do we resolve the two conundrums? In the first we have three real quantities $a=h/k$, $b=(\log n)/(\log m)$, and
$c=(c(n)-3\ep)/(c(m)+3\ep)$ and we know that $a>b\Rightarrow a>c$. From $b>a,a>c$ we would get $b>c$ by transitivity. However, in 
our situation also $a>b\Rightarrow a>c$  implies $b\ge c$, via a more subtle argument relying on the density of $\Q$ in $\R$.
The point is that we may select $a$ larger than $b$ and as close to $b$ as we wish. Assume for contrary that $c>b$. Then we select $a$ 
in-between as $c>a>b$, and $a>b\Rightarrow a>c$ gives $a>c$, a contradiction. Thus $b\ge c$. The second conundrum is more psychological
and stems from assuming $\ep>0$ to be a fixed thing. But if we drop it and regard $\ep$ as a variable on par with $m,n$, everything 
is clear. We know that $|b(m)-c(m)|,|b(n)-c(n)|<\ep\Rightarrow|\frac{c(n)}{\log n}-\frac{c(m)}{\log m}|\le3\ep(\frac{1}{\log m}+
\frac{1}{\log n})$. Thus for $m,n\in\N$ (and $m,n\ge2$) we just set $\ep=|b(m)-c(m)|+|b(n)-c(n)|$ and the implication yields the stated
conclusion (perturbing $g$ a little bit we may assume that $|b(n)-c(n)|>0$ for every $n\in\N$). 

Birch's article \cite{birc} is cited in \cite{allo,elli,elli_pnt,erdo,deko_rev,deko,deko_luca,kova,math,nark,nark_e}. 

It all started when I read the recent preprint of Shiu \cite{shiu} that reproves Segal's result \cite{sega1,sega2} that 
Euler's function $\varphi(n)$ does not have non-decreasing normal order, as a corollary of the next nice theorem.

\bigskip\noindent
{\bf Theorem (Shiu, 2016; Segal, 1964). }{\em If $f:\N\to\R_{\ge0}$ has a non-decreasing normal order, $f(n)=O(n)$, and
$\sum_{n\le x}f(n)\sim Ax^2/2$ and $\sum_{n\le x}f(n)^2\sim Bx^3/3$  as $x\to+\infty$ for some constants $A,B>0$, then
$A^2\ge B$.
}

\bigskip\noindent
For $f(n)=\varphi(n)$ (which is $O(n)$) we have $A=\prod_p(1-p^{-2})$ and $B=\prod_p(1-2p^{-2}+p^{-3})$ (see \cite{shiu} for 
proofs of these average orders). Since $A^2<B$, we conclude that $\varphi(n)$ does not have non-decreasing normal order. It follows also 
from Birch's theorem, since $\varphi(n)$ is multiplicative (and unbounded). For results on sets where $\varphi(n)$ itself is monotonous 
see Pollack, Pomerance, and Trevi\~{n}o \cite{poll}.

Finally, I was inspired by all this and the discussion at \cite{math_over} to pose the following problem.

\bigskip\noindent
{\bf Problem (MK, 2016). }{\em Does $\varphi(n)$ have an effective normal order? That is, is there a function $g:\;\N\to\N$
such that for every $\ep>0$, $\#(n\le x\;|\;\frac{\varphi(n)}{g(n)}\not\in(1-\ep,1+\ep))=o(x)$ as $x\to+\infty$, and 
$$
\mbox{one can compute $n\mapsto g(n)$ in time polynomial in $\log n$}\;?
$$
}

\medskip\noindent
{\sc Charles University, KAM MFF UK, Malostransk\'e n\'am. 25, 11800 Praha, Czechia}


\begin{thebibliography}{10}

\bibitem{allo} J.-P. Allouche, M. Mend\`es France, and  J. Peyri\`ere, Automatic Dirichlet series, J. Number Theory 81 (2000)
359--373.

\bibitem{birc} B.\,J. Birch, Multiplicative functions with non-decreasing normal order, J. London Math. Soc. 42 (1967) 149--151.

\bibitem{elli}  P.\,D.\,T.\,A. Elliott, On a conjecture of Narkiewicz about functions with non-decreasing normal order, Colloq. 
Math. 36 (1976) 289--294.

\bibitem{elli_pnt}  P.\,D.\,T.\,A. Elliott, {\em Probabilistic Number Theory. I. Mean-value Theorems}, Springer-Verlag, 
New York--Berlin, 1979.

\bibitem{elli_pr} P.\,D.\,T.\,A. Elliott, {\em Arithmetic Functions and Integer Products}, Springer-Verlag, New York, 1985.

\bibitem{erdo}  P. Erd\H os and C. Ryavec, A characterization of finitely monotonic additive functions, J. London Math. Soc. 5 
(1972), 362--367.

\bibitem{deko_rev} J.-M. de Koninck, Review of \cite{elli_pr} and \cite{mcca}, Bull. Amer. Math. Soc. 18 (1988), 230--247.

\bibitem{deko} J.-M. de Koninck, N. Doyon, and P. Letendre, On the proximity of additive and multiplicative functions, Funct. 
Approx. Comment. Math. 52 (2015) 327--344.

\bibitem{deko_luca} J.-M. de Koninck and F. Luca, {\em Analytic Number Theory. Exploring the Anatomy of Integers}, American 
Mathematical Society, Providence, RI, 2012.

\bibitem{kova} K. Kov\'acs, On the characterization of additive and multiplicative functions, Studia Sci. Math. Hungar. 18 
(1982) 1--11.

\bibitem{math} L. Matthiesen, Correlations of the divisor function, Proc. Lond. Math. Soc. 104 (2012) 827--858.

\bibitem{mcca} P.\,J. McCarthy, {\em Introduction to Arithmetical Functions}, Springer-Verlag, New York, 1986. 

\bibitem{nark} W. Narkiewicz, {\em Teoria liczb}, Pa\'nstwowe Wydawnictwo Naukowe, Warszawa, 1990 (in Polish).

\bibitem{nark_e} W. Narkiewicz, {\em Number Theory}, World Scientific Publishing Co., Singapore, 1983 (translated from the 
1977 edition of \cite{nark} by S. Kanemitsu).

\bibitem{poll}  P. Pollack, C. Pomerance, and E. Trevi\~{n}o, Sets of monotonicity for Euler's totient function, Ramanujan 
J. 30 (2013) 379--398.

\bibitem{sega1}  S.\,L. Segal, A note on normal order and the Euler $\varphi$-function, J. London Math. Soc. 39 (1964) 400--404. 

\bibitem{sega2}  S.\,L. Segal, On non-decreasing normal orders, J. London Math. Soc. 40 (1965) 459--466.

\bibitem{shiu} P. Shiu, On functions without a normal order, preprint, arXiv:1606.04533, June 2016, 4 pages. 

\bibitem{math_over} How hard is it to compute the Euler totient function?,\\ http://mathoverflow.net/questions/3274/

\end{thebibliography}
\end{document}